\input amstex
\magnification=\magstep1
\baselineskip=13pt
\documentstyle{amsppt}
\vsize=8.7truein
\NoRunningHeads

\def\conv{\operatorname{conv}}
\def\vl{\operatorname{vol}}
\def\codim{\operatorname{codim}}
\topmatter
\title Convex Geometry of Orbits \endtitle
\author Alexander Barvinok and Grigoriy Blekherman \endauthor
\address Department of Mathematics, University of Michigan, Ann Arbor,
MI 48109-1109, USA \endaddress
\email barvinok, gblekher$\@$umich.edu  \endemail
\date December 2003 \enddate
\thanks This research was partially supported by NSF Grant DMS 9734138.
\endthanks
\abstract We study metric properties of convex bodies $B$ 
and their polars $B^{\circ}$, where $B$ is the convex hull of an orbit 
under the action of a compact group $G$. Examples include the Traveling 
Salesman Polytope in polyhedral combinatorics 
($G=S_n$, the symmetric group), the set of non-negative polynomials in 
real algebraic 
geometry ($G=SO(n)$, the special orthogonal group), and 
the convex hull of the Grassmannian and the unit comass ball
in the theory of calibrated geometries ($G=SO(n)$, but with a
different action). We compute the radius of the largest ball 
contained in the symmetric Traveling Salesman Polytope, give a reasonably 
tight estimate for the radius of the Euclidean ball containing the 
unit comass ball and review (sometimes with simpler and unified proofs)   
recent results on the structure of the set of non-negative 
polynomials (the radius of the inscribed ball, 
volume estimates, and relations to the sums of squares).
Our main tool is a new simple description of the ellipsoid of the largest 
volume contained in $B^{\circ}$.  
\endabstract
\endtopmatter
\document

\head 1. Introduction and Examples \endhead

Let $G$ be a compact group acting in a finite-dimensional real vector space
$V$ and let $v \in V$ be a point. The main object of this paper is the 
convex hull
$$B=B(v)=\conv \bigl(gv: \ g \in G \bigr)$$ 
of the orbit as well as its polar 
$$B^{\circ}=B^{\circ}(v)=\bigl\{\ell \in V^{\ast}: \quad \ell(gv) \leq 1
 \quad 
\text{for all} \quad g \in G \bigr\}.$$ 
Objects such as $B$ and $B^{\circ}$ appear in many different contexts.
We give three examples below.
\example{(1.1) Example: Combinatorial optimization polytopes}
Let $G=S_n$ be the symmetric group, that is, the group of permutations 
of $\{1, \ldots, n\}$. Then $B(v)$ is a polytope and varying $V$ and 
$v$, one can obtain various polytopes of interest in combinatorial 
optimization. This idea is due to A.M. Vershik (see \cite{BV88}) and
some polytopes of this kind were studied in \cite{Ba92}.

Here we describe perhaps the most famous polytope in this family,
the Traveling Salesman Polytope 
(see, for example, Chapter 58 of \cite{Sc03}),
which exists in two major versions, symmetric and asymmetric.
Let $V$ be the space of $n \times n$ real matrices $A=(a_{ij})$ and 
let $S_n$ act in $V$ by simultaneous permutations of rows and columns:
$(ga)_{ij}=a_{g^{-1}(i) g^{-1}(j)}$ (we assume 
that $n \geq 4$). Let us choose $v$ such that 
$v_{ij}=1$ provided $|i-j|=1 \mod n$ and $v_{ij}=0$ otherwise.
Then, as $g$ ranges over the symmetric group $S_n$, matrix $gv$ ranges 
over the adjacency matrices of Hamiltonian cycles in a complete 
undirected graph 
with $n$ vertices. The convex hull $B(v)$ is called the 
{\it symmetric Traveling Salesman Polytope} (we denote it 
by $ST_n$). It has 
$(n-1)!/2$ vertices and its dimension is $(n^2-3n)/2$.

Let us choose $v \in V$ such that $v_{ij}=1$ provided $i-j=1 \mod n$   
and $v_{ij}=0$ otherwise. Then, as $g$ ranges over the symmetric group 
$S_n$, matrix $gv$ ranges over the adjacency matrices of Hamiltonian 
circuits in a complete directed graph with $n$ vertices. The convex 
hull $B(v)$ is called the {\it asymmetric Traveling Salesman Polytope}
(we denote it by $AT_n$).
It has $(n-1)!$ vertices and its dimension is $n^2-3n+1$. 
\endexample

A lot of effort has been
put into understanding of the facial structure of the 
symmetric and asymmetric Traveling Salesman Polytopes,
in particular, what are the linear inequalities that define the facets 
of $AT_n$ and $ST_n$, see Chapter 58 of \cite{Sc03}.
It follows from the computational complexity
theory that in some sense one {\it cannot} describe efficiently the 
facets of the Traveling Salesman Polytope. More precisely, if 
{\it NP} $\ne$ {\it co-NP} (as is widely believed), then there is no 
polynomial 
time algorithm, which, given an inequality, decides if it determines 
a facet of the Traveling Salesman Polytope, symmetric 
or asymmetric, see, for example, Section 5.12 of \cite{Sc03}. 
In a similar spirit,
Billera and Sarangarajan proved that any 0-1 polytope (that is,
a polytope whose vertices are 0-1 vectors), appears as a face of some 
$AT_n$ (up to an affine equivalence) \cite{BS96}. 
 
\example{(1.2) Example: Non-negative polynomials} 
Let us fix positive integers $n$ and $k$. We are interested in 
homogeneous polynomials $p: {\Bbb R}^n \longrightarrow {\Bbb R}$ of 
degree $2k$ that are non-negative for all $x=(x_1, \ldots, x_n)$. 
Such polynomials form a convex cone and we consider its compact base:
$$\aligned 
Pos_{2k,n}=\Bigl\{p: \quad &p(x) \geq 0 \quad \text{for all}
\quad x \in {\Bbb R}^n \quad \text{and} 
\\ &\int_{{\Bbb S}^{n-1}}p(x) \ dx =1 \Bigr\}, \endaligned \tag1.2.1$$
where $dx$ is the rotation invariant probability measure on the 
unit sphere ${\Bbb S}^{n-1}$.

It is not hard to see that $\dim Pos_{2k,n}={n+2k-1 \choose 2k}-1$. 

It is convenient to consider a translation $Pos_{2k,n}'$,
$p \longmapsto p-(x_1^2 + \ldots +x_n^2)^k$ of $Pos_{2k,n}$:
$$\aligned Pos_{2k,n}'=\Bigl\{p: 
\quad &p(x) \geq -1 \quad \text{for all} \quad x \in {\Bbb R}^n 
\quad \text{and} \\
&\int_{{\Bbb S}^{n-1}} p(x) \ dx=0 \Bigr\}. \endaligned \tag1.2.2$$
Let $U_{m,n}$ be the real vector space of all 
homogeneous polynomials $p: {\Bbb R}^n \longrightarrow {\Bbb R}$ 
of degree $m$ such that the average value of $p$ on ${\Bbb S}^{n-1}$ is 
0. Then, for $m=2k$, the set $Pos_{2k,n}'$ is a full-dimensional
convex body in $U_{2k,n}$.

One can view $Pos_{2k,n}'$ as the {\it negative polar} 
$-B^{\circ}(v)$ of some orbit.  

We  
consider the $m$-th tensor power 
$\left({\Bbb R}^n\right)^{\otimes m}$ of ${\Bbb R}^n$, which we view as 
the vector space of all 
$m$-dimensional arrays  
$\bigl(x_{i_1, \ldots, i_m}: 1\leq i_1, \ldots, i_m \leq n \bigr)$.
For $x \in {\Bbb R}^n$, let 
$y=x^{\otimes m}$ be the tensor with the coordinates
$y_{i_1, \ldots, i_m}=x_{i_1} \cdots x_{i_m}$.
The group $G=SO(n)$ of orientation 
preserving orthogonal transformations of ${\Bbb R}^n$ acts in 
$\left({\Bbb R}^n\right)^{\otimes m}$ by the $m$-th tensor power of its 
natural action in ${\Bbb R}^n$. In particular, 
$gy=(gx)^{\otimes m}$ for $y=x^{\otimes m}$.

Let 
us choose $e \in {\Bbb S}^{n-1}$ and let 
$w=e^{\otimes m}$. 
Then the orbit $\{gw:\  g \in G \}$ consists of the tensors 
$x^{\otimes m}$, where $x$ ranges over the unit sphere in ${\Bbb R}^n$.
The orbit $\{gw:\ g \in G\}$ lies in the symmetric part of 
$\left({\Bbb R}^n\right)^{\otimes m}$.
Let $\displaystyle q=\int_{{\Bbb S}^{n-1}} gw \ dg$ be the center 
of the orbit. If $m$ is odd then $q=0$ and if $m=2k$ is even then 
$q$ is a positive multiple of $(x_1^2 + \ldots + x_n^2)^k$.
We translate the orbit by shifting $q$ to the origin, so in 
the end we consider the convex hull $B$ of the orbit of 
$v=w-q$:
$$B=\conv\bigl(gv: \quad g \in G \bigr).$$ 

A homogeneous polynomial 
$$p(x_1, \ldots, x_n)=\sum_{1 \leq i_1, \ldots, i_m \leq n} 
c_{i_1, \ldots, i_m} x_{i_1} \cdots x_{i_m}$$
of degree $m$, viewed as a function
on the unit sphere in ${\Bbb R}^n$, is identified with
the restriction onto the orbit $\{gw:\ g \in G \bigr\}$ of  
the linear functional 
$\ell: \left({\Bbb R}^n\right)^{\otimes m} \longrightarrow {\Bbb R}$ defined
by the coefficients $c_{i_1, \ldots, i_m}$. Consequently,
the linear functionals $\ell$
on $B$ are in one-to-one correspondence with the 
polynomials $p \in U_{m,n}$. Moreover, for 
$m=2k$, the negative 
polar $-B^{\circ}$ is identified with $Pos_{2k, n}'$. If $m$ is odd, 
then $B^{\circ}=-B^{\circ}$ is the set of polynomials $p$ such 
that $|p(x)| \leq 1$ for all $x \in {\Bbb S}^{n-1}$.
\endexample

The facial structure of $Pos_{2k,n}$ is well-understood if $k=1$ or 
if $n=2$, see, for example, Section II.11 (for $n=2$) and Section II.12 
(for $k=1$) of \cite{Ba02a}. In particular, for $k=1$, the set 
$Pos_{2, n}$ is the convex body of positive semidefinite 
$n$-variate quadratic forms of trace $n$. The faces of 
$Pos_{2,n}$ are parameterized by the subspaces of ${\Bbb R}^n$: 
if $L \subset {\Bbb R}^n$ is a subspace then the corresponding face 
is 
$$F_L=\Bigl\{p \in Pos_{2,n}: \quad p(x)=0 \quad \text{for all} \quad 
x \in L \Bigr\}$$
and $\dim F_L=r(r+1)/2-1$, where $r=\codim L$. Interestingly, 
for large $n$, the set $Pos_{2,n}$ is a counterexample to famous Borsuk's
conjecture \cite{K95}.

For any $k \geq 2$, the situation is much more complicated:
 the {\it membership problem}
for $Pos_{2k,n}$: 
\medskip
given a polynomial, decide whether it belongs to 
$Pos_{2k, n}$, 
\medskip
\noindent is NP-hard, which indicates that the facial structure of 
$Pos_{2k,n}$ is probably hard to describe.

\example{(1.3) Example: Convex hulls of Grassmannians and 
calibrations} Let $G_m({\Bbb R}^n)$ be the Grassmannian of all 
oriented $m$-dimensional subspaces of ${\Bbb R}^n$. Let us consider 
$G_m({\Bbb R}^n)$ as a subset of $V_{m,n}=\bigwedge^m {\Bbb R}^n$
via the Pl\"ucker embedding. Namely, let $e_1, \ldots, e_n$ be the 
standard basis of ${\Bbb R}^n$. We make $V_{m,n}$ a 
Euclidean space by choosing an orthonormal basis
$e_{i_1} \wedge \ldots \wedge e_{i_m}$ for 
$1 \leq i_1 < \ldots < i_m \leq n$.
 Thus the coordinates of a subspace 
$x \in G_m({\Bbb R}^n)$ are indexed by $m$-subsets 
$1 \leq i_1 < i_2 < \ldots < i_m \leq n$ of $\{1, \ldots, n\}$  
and the coordinate $x_{i_1, \ldots, i_m}$ is equal to the oriented volume
of the parallelepiped spanned by the orthogonal projection of 
$e_{1_1}, \ldots, e_{i_m}$ onto $x$. This identifies $G_m({\Bbb R}^n)$ 
with 
a subset of the unit sphere in $V_{m,n}$. 
The convex hull $B=\conv\left(G_m({\Bbb R}^n)\right)$, called 
the {\it unit mass ball}, turns out to be of interest in the 
theory of calibrations and area-minimizing surfaces: a face of 
$B$ gives rise to a family of 
$m$-dimensional
area-minimizing surfaces whose tangent planes belong to the face,
see \cite{HL82} and \cite{F88}. The {\it comass} of 
a linear functional $\ell: V_{m,n} \longrightarrow {\Bbb R}$ is the
maximum value of $\ell$ on $G_m({\Bbb R}^n)$.
A {\it calibration} is a linear functional
$\ell: V_{m,n} \longrightarrow {\Bbb R}$ of comass 1. 
The polar $B^{\circ}$ is called the {\it unit comass ball}.

One can easily view 
$G_m({\Bbb R}^n)$ as an orbit. We let $G=SO(n)$, the group of
orientation-preserving orthogonal transformations of ${\Bbb R}^n$,
 and consider the action of 
$SO(n)$ in $V_{m,n}$ by the $m$-th exterior power of its defining action
in ${\Bbb R}^n$. Choosing $v=e_1 \wedge \ldots \wedge e_m$, 
we observe that $G_m({\Bbb R}^n)$ is the orbit 
$\{gv: \ g \in G \}$. It is easy to see that 
$\dim \conv\bigl(G_m({\Bbb R}^n)\bigr)={n \choose m}$.

This example was suggested to the authors by B. Sturmfels and J. Sullivan.
\endexample

The facial structure of the convex hull of $G_m({\Bbb R}^n)$ is understood 
for $m \leq 2$, for $m \geq n-2$ and for some special values of $m$ and $n$,
see \cite{HL82}, \cite{HM86} and \cite{F88}. If $m=2$, then 
the faces of the unit mass ball are as follows: let us choose an 
even-dimensional subspace $U \subset {\Bbb R}^m$ and an orthogonal 
complex structure on $U$, thus identifying $U={\Bbb C}^{2k}$ for 
some $k$. Then the corresponding face of 
$\conv \left( G_m ({\Bbb R}^n)\right)$ is the convex 
hull of all oriented planes 
in $U$ identified with complex lines in ${\Bbb C}^{2k}$. 

In general, it appears to be difficult to
describe the facial structure of the unit mass ball.
The authors do not know the complexity status of the 
{\it membership problem} for the unit mass ball: 
\medskip
given a point $x \in \bigwedge^m {\Bbb R}^n$, decide if it lies in 
$\conv\left(G_m({\Bbb R}^n)\right)$,
\medskip
but suspect that the problem is NP-hard if $m \geq 3$ is fixed and 
$n$ is allowed to grow. 
\bigskip
The above examples suggest that the boundary of $B$ and $B^{\circ}$ can 
get very complicated, so there is little hope in understanding the 
combinatorics (the facial structure) of general convex hulls of orbits and 
their polars. Instead, we study metric properties of convex hulls.
Our approach is through approximation of a complicated convex body by 
a simpler one. 

As is known, every convex body contains a unique ellipsoid 
$E_{\max}$ of the maximum volume and is contained in a unique ellipsoid 
$E_{\min}$ of the minimum volume, see \cite{B97}. Thus ellipsoids 
$E_{\max}$ and $E_{\min}$ provide reasonable ``first approximations''
to a convex body. 

The main result of Section 2 is Theorem 2.4
which states that the maximum volume ellipsoid of $B^{\circ}$ 
consists of the linear functionals $\ell: V \longrightarrow {\Bbb R}$ 
such that the average value of $\ell^2$ on the orbit does 
not exceed $(\dim V)^{-1}$.  We compute the minimum- and 
maximum- volume ellipsoids of the symmetric Traveling 
Salesman Polytope, which both turn out 
to be balls under the ``natural'' Euclidean metric and 
ellipsoid $E_{\min}$ of the asymmetric Traveling Salesman 
Polytope, which turns out to be slightly stretched in the 
direction of the skew-symmetric matrices. As an immediate corollary of
Theorem 2.4, we obtain the description of the maximum volume 
ellipsoid of the set of non-negative polynomials (Example 1.2), as
a ball of radius $\left( {n+2k-1 \choose 2k} -1 \right)^{-1/2}$ in the 
$L^2$-metric. We also compute the minimum volume 
ellipsoid of the convex hull of the Grassmannian and hence the 
maximum volume ellipsoid of the unit comass ball (Example 1.3).

In Section 3, we obtain some inequalities which allow us to 
approximate the maximum value of a linear functional $\ell$ 
on the orbit by an $L^p$-norm of $\ell$. We apply those inequalities
in Section 4. We obtain a reasonably tight estimate of 
the radius of the Euclidean ball containing the unit comass ball and
show that the classical K\"ahler and special Lagrangian faces of 
the Grassmannian, are, in fact, rather ``shallow'' (Example 1.3).
Also, we review (with some proofs and some sketches) the recent results of 
\cite{Bl03}, which show that for most values of $n$ and $k$ 
the set of non-negative $n$-variate polynomials
of degree $2k$ is much larger than its subset consisting of the sums 
of squares of polynomials of degree $k$.     

\head 2. Approximation by Ellipsoids \endhead

Let $B \subset V$ be a convex
body in a finite-dimensional real vector space.
We assume that $\dim B=\dim V$. Among all ellipsoids contained 
in $B$ there is a unique ellipsoid $E_{\max}$ of the maximum volume,
which we call the maximum volume ellipsoid of $B$ and which is also called
the John ellipsoid of $B$ or the L\"owner-John ellipsoid 
of $B$. Similarly, among all ellipsoids containing 
$B$ there is a unique ellipsoid $E_{\min}$ of the minimum volume, which 
we call the minimum volume ellipsoid of $B$ and which is also called the 
L\"owner or the L\"owner-John ellipsoid. The maximum and minimum volume 
ellipsoids of $B$ do not depend on the volume form chosen in $V$, they 
are intrinsic to $B$.

Assuming that the center of $E_{\max}$ is the origin, we have 
$$E_{\max} \subset B \subset \left(\dim B\right) E_{\max}.$$
If $B$ is symmetric about the origin, that is, if $B=-B$ then 
the bound can be strengthened:
$$E_{\max} \subset B \subset \left(\sqrt{\dim B}\right) E_{\max}.$$
More generally, let us suppose that $E_{\max}$ is centered at the origin.
The {\it symmetry coefficient} of $B$ with respect 
to the origin is the largest $\alpha>0$ such that
$-\alpha B \subset B$. Then we have
$$E_{\max} \subset B \subset \left(\sqrt{\dim B \over \alpha}\right)
 E_{\max},$$
where $\alpha$ is the symmetry coefficient of $B$ with respect to the 
origin.

Similarly, assuming that $E_{\min}$ is centered at the origin, 
we have 
$$\left(\dim B\right)^{-1} E_{\min} \subset B \subset E_{\min}.$$
If, additionally, $\alpha$ is the symmetry coefficient of $B$ with respect
to the origin, then
$$\left(\sqrt{\alpha \over \dim B}\right) 
E_{\min} \subset B \subset E_{\min}.$$
In particular, if $B$ is symmetric about the origin, then 
$$\left(\dim B\right)^{-1/2} E_{\min} \subset B \subset E_{\min}.$$
These, and other interesting properties of the minimum- and maximum-
volume ellipsoids can be found in \cite{B97}, see also the original 
paper \cite{J48}, \cite{Bl03}, and Chapter V of \cite{Ba02b}.     

Suppose that a compact group $G$ acts in $V$ by linear transformations and
that $B$ is invariant under the action: $gB=B$ for all $g \in G$.
Let $\langle \cdot, \cdot \rangle$ be a $G$-invariant scalar product 
in $V$, so $G$ acts in $V$ by isometries.
Since the ellipsoids $E_{\max}$ and $E_{\min}$ associated with $B$ 
are unique, they also have to be invariant under the action of $G$.
If the group of symmetries of $B$ is sufficiently rich, we may be able to 
describe $E_{\max}$ or $E_{\min}$ precisely. 

The following simple observation will be used throughout this section.
Let us suppose that the action of $G$ in $V$ is irreducible:
if $W \subset V$ is a $G$-invariant subspace, then either 
$W=\{0\}$ or $W=V$. Then, the ellipsoids $E_{\max}$ and $E_{\min}$ of 
a $G$-invariant convex body $B$ are necessarily balls centered at the origin:
$$E_{\max}=\bigl\{x \in V:\ \langle x, x \rangle \leq r^2 \bigr\}
\quad \text{and} \quad 
E_{\min}=\bigl\{x \in V: \ \langle x, x \rangle \leq R^2 \bigr\}$$
for some $r, R>0$. 

Indeed, since the action of $G$ is irreducible, the origin is the 
only $G$-invariant point and hence both $E_{\max}$ and $E_{\min}$ 
must be centered at the origin. Furthermore, an ellipsoid $E \subset V$ 
centered at the origin is defined by the inequality 
$E=\bigl\{x: q(x) \leq 1 \bigr\}$, where $q: V \longrightarrow {\Bbb R}$ 
is a positive definite quadratic form. If $E$ is $G$-invariant, then 
$q(gx)=q(x)$ for all $g \in G$ and hence the eigenspaces of $q$ must 
be $G$-invariant. Since the action of $G$ is irreducible, there is 
only one eigenspace which coincides with $V$, from which 
$q(x)=\lambda \langle x,x \rangle$ for some $\lambda >0$ and all $x \in V$
and $E$ is a ball.  
 
This simple observation allows us to compute 
ellipsoids $E_{\max}$ and 
$E_{\min}$ of the Symmetric Traveling Salesman Polytope (Example 1.1).

\example{(2.1) Example: The minimum and maximum volume ellipsoids of the 
symmetric Traveling Salesman Polytope}
In this case, $V$ is the space of $n \times n$ real matrices, on which 
the symmetric group $S_n$ acts by simultaneous permutations of rows and
columns, see Example 1.1. Let us introduce an $S_n$-invariant scalar 
product by 
$$\big\langle a, b \big\rangle=\sum_{i,j=1}^n a_{ij} b_{ij} \quad 
\text{for} \quad a=(a_{ij}) \quad \text{and} \quad b=(b_{ij})$$
and the corresponding Euclidean norm $\|a\|=\sqrt{\langle a, a \rangle}$.
It is not hard to see that the affine hull of the 
symmetric Traveling Salesman Polytope $ST_n$ consists
of the symmetric matrices with 0 diagonal and row and column sums equal to
2, from which one can deduce the formula $\dim ST_n=(n^2-3n)/2$.
Let us make the affine hull of $ST_n$ a vector space by choosing the origin 
at $c=(c_{ij})$ with $c_{ij}=2/(n-1)$ for $i \ne j$ and $c_{ii}=0$,
the only fixed point of the action. 
One can see that the action of $S_n$ on the affine hull of $ST_n$
is irreducible and corresponds to the Young diagram $(n-2,2)$, see,
for example, Chapter 4 of \cite{FH91}.

Hence the maximum- and minimum- volume ellipsoids 
of $ST_n$ must be balls in the affine hull of $ST_n$ centered at 
$c$. Moreover, since the boundary of the minimum volume ellipsoid 
$E_{\min}$ must contain the vertices of $ST_n$, we conclude that the 
radius 
of the ball representing $E_{\min}$ is equal to $\sqrt{2n(n-3)/(n-1)}$.

One can compute the symmetry coefficient of $ST_n$ with respect to 
the center $c$. Suppose that $n \geq 5$.
 Let us choose a vertex $v$ of $ST_n$ and let us 
consider the functional $\ell(x)=\langle v-c, x-c \rangle$ on 
$ST_n$. The maximum value of $2n(n-3)/(n-1)$ is attained at 
$x=v$ while the minimum value of $-4n/(n-1)$ is attained at 
the face $F_v$ of $ST_n$ with the vertices $h$ such that 
$\langle v, h \rangle=0$ (combinatorially, $h$ correspond to 
Hamiltonian cycles in the graph obtained from the complete 
graph on $n$ vertices by deleting the edges of the Hamiltonian 
cycle encoded by $v$). Moreover, one can show that for 
$\lambda=2/(n-3)$, we have $-\lambda(v-c)+c \in F_v$. 
This implies that the coefficient of symmetry of $ST_n$ with 
respect to $c$ is equal to $2/(n-3)$. Therefore $ST_n$ contains 
the ball centered at $c$ and of the radius 
$\sqrt{8/\bigl((n-1)(n-3)\bigr)}$ (for $n \geq 5$).

One can observe that the ball centered at $c$ and of the radius 
$\sqrt{8/\bigl((n-1)(n-3)\bigr)}$ touches the boundary of $ST_n$.
Indeed, let $b=(b_{ij})$ be the centroid of the set of vertices
$x$ of $ST_n$ with $x_{12}=x_{21}=0$. Then 
$$b_{ij}=\cases 0 &\text{if\ } 1 \leq i, j \leq 2 \\
{2 \over n-2} &\text{if\ } i=1,2 \ \text{and} \ j >2 \quad \text{or} 
\quad j=1,2 \ \text{and} \ i >2 \\
2(n-4) \over {(n-2)(n-3)} &\text{if\ } i,j\geq 3, \endcases$$
and the distance from $c$ to $b$ 
is precisely $\sqrt{8/\bigl((n-1)(n-3)\bigr)}$.

Hence for $n \geq 5$ 
the maximum volume ellipsoid $E_{\max}$ is the ball centered 
at $c$ of the radius $\sqrt{8/\bigl((n-1)(n-3)\bigr)}$. 

Some bounds on the radius of the largest inscribed ball for 
a polytope from a particular family of 
combinatorially defined polytopes are computed in \cite{V95}. 
The family of polytopes includes the symmetric Traveling Salesman Polytope,
although in its case the bound from \cite{V95} is not optimal.  
\endexample

If the action of $G$ in the ambient space $V$ is not irreducible, 
the situation is more complicated. For one thing, there  
is more than one (up to a scaling factor) $G$-invariant scalar 
product, hence the notion of a ``ball'' is not really defined.
However, we are still able to describe the minimum volume ellipsoid 
of the convex hull of an orbit.

Without loss of generality, we assume that the orbit 
$\bigl\{gv:\ g \in G\bigr\}$ spans $V$ affinely.
Let $\langle \cdot, \cdot \rangle$ be a $G$-invariant scalar product in 
$V$. As is known, $V$ can be decomposed into the direct sum of 
pairwise orthogonal invariant subspaces $V_i$, such that 
the action of $G$ in each $V_i$ is irreducible. It is important 
to note that the decomposition is {\it not} unique: non-uniqueness 
appears when some of $V_i$ are isomorphic, that, is, when there exists 
an isomorphism $V_i \longrightarrow V_j$ which commutes with $G$.  
If the decomposition is unique, we say that the action of $G$ is 
{\it multiplicity-free}.

Since the orbit spans $V$ 
affinely, the orthogonal projection $v_i$ of $v$ onto each 
$V_i$ must be non-zero (if $v_i=0$ then the orbit lies in 
$V_i^{\bot}$). Also, the origin in $V$ must be 
the only invariant point of the action of $G$ (otherwise, the orbit is 
contained in the hyperplane $\langle x, u \rangle=\langle v, u \rangle$, 
where $u \in V$ is a non-zero vector fixed by the action of $G$).

\proclaim{(2.2) Theorem} Let $B$ be the convex hull of the orbit 
of a vector $v \in V$:
$$B=\conv\Bigl(gv:\ g \in G \Bigr).$$ 
Suppose that the affine hull of $B$ is $V$.
  
Then there exists a decomposition
$$V=\bigoplus_{i} V_i$$
of $V$ into the direct sum of pairwise orthogonal 
irreducible components with the following properties.
  
The minimum volume ellipsoid 
$E_{\min}$ of $B$ is defined by the inequality
$$
E_{\min}=\Bigl\{x: \quad \sum_i {\dim V_i \over \dim V} \cdot
{\langle x_i, x_i \rangle \over \langle v_i, v_i \rangle} \leq 1
\Bigr\},\tag2.2.1$$
where $x_i$ (resp. $v_i$) is the orthogonal projection of $x$ 
(resp. $v$) onto $V_i$. 

We have
$$\int_G \langle x, gv \rangle^2 \ dg=\sum_i {\langle x_i, x_i \rangle
\langle v_i, v_i \rangle \over \dim V_i} \quad \text{for all} \quad 
x \in V, \tag2.2.2$$
where $dg$ is the Haar probability measure on $G$.
\endproclaim
\demo{Proof} Let us consider the quadratic form 
$q: V \longrightarrow {\Bbb R}$ defined by 
$$q(x)=\int_G \langle x, gv \rangle^2 \ dg.$$
We observe that $q$ is $G$-invariant, that is, $q(gx)=q(x)$ for all 
$x \in V$ and all $g \in G$. Therefore, the eigenspaces of $q$ are
$G$-invariant. Writing the eigenspaces as direct sums of pairwise 
orthogonal invariant subspaces where the action of $G$ is irreducible,
we obtain a decomposition $V=\bigoplus_i V_i$ such that 
$$q(x)=\sum_i \lambda_i \langle x_i, x_i \rangle \quad 
\text{for all} \quad x \in V$$
and some $\lambda_i \geq 0$. Recall that $v_i \ne 0$ for all $i$ since
the orbit $\{gv: g \in G\}$ spans $V$ affinely. 

To compute $\lambda_i$, we substitute 
$x \in V_i$ and observe that the trace of 
$$q_i(x)=\int_G \langle x, gv_i \rangle^2 \ dg$$
as a quadratic form $q_i: V_i \longrightarrow {\Bbb R}$ 
is equal to $\langle v_i, v_i \rangle$. Hence we must have 
$\lambda_i=\langle v_i, v_i \rangle/\dim V_i$, which proves (2.2.2), cf.
\cite{Ba02a}.

We will also use the polarized form of (2.2.2):
$$\int_G \langle x, gv \rangle \langle y, gv \rangle \ dg =
\sum_i {\langle x_i, y_i \rangle \langle v_i, v_i \rangle \over \dim V_i},
\tag2.2.3$$
obtained by applying (2.2.2) to $q(x+y)-q(x)-q(y)$.

Next, we observe that the ellipsoid $E$ defined by the inequality 
(2.2.1) contains the orbit $\bigl\{gv:\ g \in G\bigr\}$ on 
its boundary and hence 
contains $B$. 

Our goal is to show that $E$ is the minimum volume ellipsoid.
It is convenient to introduce a new scalar product:
$$(a,b)=\sum_i {\dim V_i \over \dim V} \cdot
{\langle a_i, b_i \rangle \over \langle v_i, v_i \rangle} \quad 
\text{for all} \quad a,b \in V.$$
Obviously $(\cdot ,\cdot)$ is a $G$-invariant scalar product.
Furthermore, the ellipsoid $E$ defined by (2.2.1) is the unit ball in 
the scalar product $(\cdot, \cdot)$.

Now,
$$(c, gv)=\sum_i {\dim V_i \over \dim V} \cdot 
{\langle c_i, gv \rangle \over \langle v_i, v_i \rangle}$$
and hence
$$(c, gv)^2=\sum_{i,j} {(\dim V_i)(\dim V_j) \over (\dim V)^2} \cdot  
{\langle c_i, gv \rangle \langle c_j, gv \rangle \over 
\langle v_i, v_i \rangle^2}.$$
Integrating and using (2.2.3), we get 
$$\int_G (c, gv)^2 \ dg={1 \over \dim V}\sum_i {\dim V_i \over \dim V}
\cdot {\langle c_i, c_i \rangle \over \langle v_i, v_i \rangle}=
{(c,c)\over \dim V}. \tag2.2.4$$

Since the origin is the only 
fixed point of the action of $G$, the minimum volume ellipsoid should 
be centered at the origin. 

Let $e_1, \ldots, e_k$ for $k=\dim V$ be an orthonormal basis with respect
to the scalar product $(\cdot, \cdot)$. Suppose that $E' \subset V$ 
is an ellipsoid defined by 
$$E'=\Bigl\{x \in V:\quad \sum_{j=1}^k {(x, e_j)^2 \over \alpha_j^2}
\leq 1 \Bigr\}$$
for some $\alpha_1, \ldots, \alpha_k>0$. 
To show that $E$ is the minimum volume ellipsoid, it suffices to show 
that as long as $E'$ contains the orbit $\bigl\{gv:\ g \in G\bigr\}$, 
we must have 
$\vl E' \geq \vl E$, which is equivalent to
$\alpha_1 \cdots \alpha_k \geq 1$. 

Indeed, since $gv \in E'$, we must have
$$\sum_{j=1}^k {(e_j, gv)^2 \over \alpha_j^2} \leq 1 \quad 
\text{for all} \quad g \in G.$$
Integrating, we obtain
$$\sum_{j=1}^k {1 \over \alpha_j^2}\int_G (e_j, gv)^2 \ dg \leq 1.$$
Applying (2.2.4), we get
$${1 \over \dim V} \sum_{j=1}^k {1 \over \alpha_j^2} \leq 1.$$
Since $k=\dim V$, from the inequality between the arithmetic and 
geometric means, we get that $\alpha_1 \ldots \alpha_k \geq 1$,
which completes the proof. 
{\hfill \hfill \hfill} \qed
\enddemo
\remark{Remark} We note that in the part of the proof
where we compare the volumes of $E'$ and $E$, we reproduce 
the ``sufficiency'' (that is, ``the easy'') part of John's criterion for 
optimality of an ellipsoid, cf., for example, \cite{B97}. 
\endremark
\medskip
Theorem 2.2 allows us to compute the minimum volume ellipsoid of 
the asymmetric Traveling Salesman Polytope, see Example 1.1.
\example{(2.3) Example: the minimum volume ellipsoid of the 
asymmetric Traveling Salesman Polytope} In this case (see 
Examples 1.1 and 2.1), $V$ is the space of $n \times n$ matrices
with the scalar product and the action of the symmetric group 
$S_n$ defined as in Example 2.1.
On can observe that the affine hull of $AT_n$ consists of the matrices 
with zero diagonal and row and column sums equal to 1, from which one 
can deduce the formula $\dim AT_n=n^2-3n+1$.

The affine hull of $AT_n$ is $S_n$-invariant. 
We make the affine hull of $AT_n$ a vector space 
by choosing the origin at $c=(c_{ij})$ with $c_{ij}=1/(n-1)$ for 
$i \ne j$ and $c_{ii}=0$, the only fixed point of the action. 
The action of $S_n$ on the affine hull of $AT_n$ is reducible
and multiplicity-free, so there is no ambiguity in choosing the 
irreducible components. The affine hull is 
the sum of two irreducible invariant subspaces $V_s$ and $V_a$.

Subspace $V_s$ consists of the matrices $x+c$, where 
$x$ is a symmetric matrix with zero diagonal and 
zero row and column sums. One can see that the action of 
$S_n$ in $V_s$ is irreducible and corresponds to the Young diagram 
$(n-2,2)$, see, for example, Chapter 4 of \cite{FH91}. We have 
$\dim V_s=(n^2-3n)/2$

Subspace $V_s$ consists of the matrices matrices $x+c$, 
where $x$ is a skew-symmetric matrix with zero row and column sums.
One can see that the action of 
$S_n$ in $V_s$ is irreducible and corresponds to the Young diagram 
$(n-2,1,1)$, see, for example, Chapter 4 of \cite{FH91}. We have 
$\dim V_s=(n-1)(n-2)/2$.

The orthogonal projection onto $V_s$ is defined by $x \longmapsto (x+x^t)/2$,
while the orthogonal projection onto $V_a$ is 
defined by $x \longmapsto (x-x^t)/2 + c$.
 
Applying Theorem 2.2, we conclude that the minimum volume ellipsoid of 
$AT_n$ is defined in the affine hull of $AT_n$ by the inequality:
$$\split 
&(n-1) \sum_{1 \leq i \ne j \leq n} \Bigl({x_{ij}+x_{ji} \over 2}-
{1 \over n-1} \Bigr)^2\\ +
&{(n-1)(n-2) \over n} \sum_{1 \leq i \ne j \leq n}
\Bigl({x_{ij}-x_{ji} \over 2} \Bigr)^2 \leq n^2-3n+1.\endsplit$$
Thus one can say that the minimum volume ellipsoid of the asymmetric 
Traveling Salesman Polytope is slightly stretched in the direction 
of skew-symmetric matrices. 
\endexample

The dual version of Theorem 2.2 is especially simple.
\proclaim{(2.4) Theorem} Let $G$ be a compact group acting in a
finite-dimensional real vector space $V$. Let $B$ be the convex hull of 
the orbit of a vector $v \in V$:
$$B=\conv\Bigl(gv: \ g \in G\Bigr).$$
Suppose that the affine hull of $B$ is $V$.

Let $V^{\ast}$ be the dual to $V$ and let 
$$B^{\circ}=\Bigl\{\ell \in V^{\ast}: \quad \ell(x) \leq 1 
\quad \text{for all}
\quad x\in B \Bigr\}$$
be the polar of $B$. Then the maximum volume ellipsoid of $B^{\circ}$ is 
defined by the inequality
$$E_{\max}=\Bigl\{\ell \in V^{\ast}: \quad \int_G \ell^2(gv) \ dg \leq 
{1 \over \dim V} \Bigr\}.$$  
\endproclaim
\demo{Proof} Let us introduce a $G$-invariant scalar product
$\langle \cdot, \cdot \rangle$ in $V$, thus identifying $V$ and 
$V^{\ast}$. Then 
$$B^{\circ}=\Bigl\{c \in V:\quad \langle c, gv \rangle \leq 1 
\quad \text{for all} \quad g \in G \Bigr\}.$$
Since the origin is the only point fixed by the action of $G$, the 
maximum volume ellipsoid $E_{\max}$ of $B^{\circ}$ is centered at the origin.
Therefore, $E_{\max}$ must be the polar of 
the minimum volume ellipsoid of $B$. 

Let $\displaystyle V=\bigoplus_i V_i$ be the decomposition of Theorem 2.2. 
Since $E_{\max}$ is the polar of the ellipsoid $E_{\min}$ associated 
with $B$, from (2.2.1), we get
$$E_{\max}=\Bigl\{c: \quad \dim V \sum_i 
{\langle c_i, c_i \rangle \langle v_i, v_i \rangle \over \dim V_i} 
\leq 1 \Bigr\}.$$
Applying (2.2.2), we get
$$E_{\max}=\Bigl\{c: \quad \int_G \langle c, gv \rangle^2 \ dg
\leq {1 \over \dim V} \Bigr\},$$
which completes the proof.
{\hfill \hfill \hfill} \qed
\enddemo

\remark{Remark} Let $G$ be a compact group acting in a finite-dimensional
real vector space $V$ and let $v \in V$ be a point such that 
the orbit $\bigl\{gv: \ g \in V \bigr\}$ spans $V$ affinely. 
Then the dual space $V^{\ast}$ acquires a natural scalar product 
$$\langle \ell_1, \ell_2 \rangle=\int_G \ell_1(gv) \ell_2(gv) \ dg$$
induced by the scalar product in $L^2(G)$.
Theorem 2.4 states that the maximum volume ellipsoid of the polar of 
the orbit is the ball of radius $(\dim V)^{-1/2}$ in this scalar product. 

By duality, $V$ acquires the dual scalar product (which we denote below by
$\langle \cdot, \cdot \rangle$ as well). It is a 
constant multiple of the 
product $( \cdot, \cdot)$ introduced in the proof of Theorem 2.2:
$\langle u_1, u_2 \rangle=(\dim V) (u_1, u_2)$. 
We have $\langle v, v \rangle=\dim V$ and the minimum volume 
ellipsoid of the convex hull of the orbit of $v$ is the ball of radius 
$\sqrt{\dim V}$. 
\endremark
\medskip
As an immediate application of Theorem 2.4, we compute the 
maximum volume ellipsoid of the set of non-negative polynomials, see
Example 1.2.
\example{(2.5) Example: the maximum volume ellipsoid of the 
set of non-negative polynomials} 
In this case, $U_{2k,n}^{\ast}$ is the space of all
homogeneous polynomials 
$p: {\Bbb R}^n \longrightarrow {\Bbb R}$ of degree $2k$ with the zero 
average 
on the unit sphere ${\Bbb S}^{n-1}$, so 
$\dim U_{2k,n}^{\ast}={n+2k-1 \choose 2k}-1$.
We view such a polynomial $p$ as
a linear functional $\ell$ on an orbit $\bigl\{gv: \ g \in G\bigr\}$ 
in the action of the orthogonal group 
$G=SO(n)$ in $\left({\Bbb R}^n\right)^{\otimes 2k}$ 
and the shifted set $Pos_{2k, n}'$ of non-negative polynomials as 
the negative polar $-B^{\circ}$ of the orbit, see Example 1.2.
In particular, under this identification $p \longleftrightarrow \ell$,
we have
$$\int_{{\Bbb S}^{n-1}} p^2(x) \ dx = \int_{G} \ell^2 (gv) \ dg,$$ 
where $dx$ and $dg$ are the Haar probability measures on 
${\Bbb S}^{n-1}$ and $SO(n)$ respectively.

Applying Theorem 2.4 to $-B^{\circ}$, we conclude that the maximum 
volume ellipsoid of $-B^{\circ}=Pos_{2k,n}'$ consists of the polynomials $p$ such
that
$$ \int_{{\Bbb S}^{n-1}} p(x) \ dx =0 \quad \text{and} \quad 
\int_{{\Bbb S}^{n-1}} p^2(x) \ dx  \leq 
\left({n+2k-1 \choose 2k}-1 \right)^{-1}.$$
Consequently, the maximum volume ellipsoid of $Pos_{2k,n}$ consists 
of the polynomials $p$ such that
$$\int_{{\Bbb S}^{n-1}} p(x) \ dx=1 
\quad \text{and} \quad \int_{{\Bbb S}^{n-1}} \left(p(x)-1\right)^2 \ dx \leq 
\left({n+2k-1 \choose 2k}-1 \right)^{-1}.$$

Geometrically, the maximum volume ellipsoid of $Pos_{2k,n}$ can 
be described as follows. Let us introduce a scalar product in the 
space of polynomials by
$$\langle f,\ g \rangle =\int_{{\Bbb S}^{n-1}} f(x) g(x) \ dx,$$
where $dx$ is the rotation invariant probability measure, as above.
Then the maximum volume ellipsoid of $Pos_{2k,n}$ is the ball 
centered at $r(x)=(x_1^2 + \ldots + x_n^2)^k$ and of the 
radius $\displaystyle \left({n+2k-1 \choose 2k}-1 \right)^{-1/2}$. 
This result was first obtained
by more direct and complicated computations in \cite{Bl02}.
In the same paper, 
G. Blekherman also determined the coefficient of symmetry of 
$Pos_{2k,n}$ (with respect to the center $r$), it turns out to be equal to 
$\left({n+k-1 \choose k}-1\right)^{-1}$.

It follows then that $Pos_{2k,n}$ is contained in the ball centered 
at $r$ and of the radius 
$\left({n+k-1 \choose k} -1 \right)^{1/2}$. This estimate is poor 
if $k$ is fixed and $n$ is allowed to grow: as follows from results of 
Duoandikoetxea \cite{D87}, 
for any fixed $k$, the set $Pos_{2k, n}$ is contained in 
a ball of a fixed radius, as $n$ grows. However, the estimate gives 
the right logarithmic order if $k \gg n$, which one can observe by 
inspecting a polynomial $p \in Pos_{2k,n}$ that is the $2k$-th power 
of a linear function.
\endexample

We conclude this section by computing the the minimum volume 
ellipsoid of the convex hull of the Grassmannian and, consequently, 
the maximum volume ellipsoid of the unit comass ball, see Example 1.3.

\example{(2.6) Example: the minimum volume ellipsoid of the convex 
hull of the Grassmannian} In this case, $V_{m,n}=\bigwedge^m {\Bbb R}^n$
with the orthonormal basis 
$e_I=e_{i_1} \wedge \ldots \wedge e_{i_m}$, where 
$I$ is an $m$-subset $1 \leq i_1 < i_2 < \ldots <i_m \leq n$
of the set $\{1, \ldots, n\}$ and $e_1, \ldots, e_n$ is the 
standard orthonormal basis of ${\Bbb R}^n$.
 
Let $\langle \cdot, \cdot \rangle$ be the corresponding scalar product
in $V_{m,n}$,
so that 
$$\langle a, b \rangle =\sum_I 
a_I b_I,$$
where $I$ ranges over all $m$-subsets of $\{1, \ldots, n\}$. 
The scalar product allows us to identify $V^{\ast}_{m,n}$ with 
$V_{m,n}$.
First, we find the maximum volume ellipsoid of the unit comass ball
$B^{\circ}$, that 
is the polar of the convex hull
$B=\conv\left(G_m({\Bbb R}^n)\right)$ of the Grassmannian. 

A linear functional $a \in V_{m,n}^{\ast}=V_{m,n}$ is 
defined by its coefficients $a_I$. To apply Theorem 2.4, we have 
to compute 
$$\int_{SO(n)} \langle a, gv \rangle^2 \ dg=
\int_{G_m({\Bbb R}^n)} \langle a, x \rangle^2 \ dx,$$
where $dx$ is the Haar probability measure on the Grassmannian 
$G_m({\Bbb R}^n)$.
We note that 
$$\int_{G_m({\Bbb R}^n)} \langle e_I, x \rangle 
\langle e_J, x \rangle \ dx=0$$
for $I \ne J$, since for $i \in I \setminus J$, the reflection 
$e_i \longmapsto -e_i$ of ${\Bbb R}^n$ induces an isometry of 
$V_{m,n}$, which maps $G_m({\Bbb R}^n)$ onto itself, reverses the 
sign of $\langle e_I, x \rangle$ and does not change 
$\langle e_J, x \rangle$. Also,
$$\int_{G_m({\Bbb R}^n)} \langle e_I, x \rangle^2\ dx={n \choose m}^{-1},$$
since the integral does not depend on $I$ and 
$\sum_I \langle e_I, x \rangle^2=1$ for all $x \in G_m({\Bbb R}^n)$.

By Theorem 2.4, we conclude that the maximum volume ellipsoid of 
the unit comass ball $B^{\circ}$ is defined by the inequality
$$E_{\max}=\Bigl\{a \in V_{m,n}: \quad \sum_I a_I^2 \leq 1 \Bigr\},$$
that is, the unit ball in the Euclidean metric of $V_{m,n}$. Since
 $B^{\circ}$ is centrally symmetric, we conclude that 
$B^{\circ}$ is contained in the ball of radius 
$\displaystyle {n \choose m}^{1/2}$. As follows from Theorem 4.1, 
this estimate is optimal up to a factor of $\bigl(mn\ln(m+1)\bigr)^{1/2}$.

Consequently, the convex hull $B$ of the Grassmannian is 
contained in the unit ball of $V_{m,n}$, which is the minimum volume 
ellipsoid of $B$, and contains a ball of radius 
$\displaystyle {n \choose m}^{-1/2}$. 
Again, the estimate of the radius of the inner 
ball is optimal up to a factor of $\bigl(mn\ln(m+1)\bigr)^{1/2}$.
\endexample 

\head 3. Higher Order Estimates \endhead

The following construction can be used to get a better understanding
of metric properties of an orbit $\bigl\{ gv: \ g \in G \bigr\}$. 
Let us choose a positive integer $k$ and let us consider the 
$k$-th tensor power 
$$V^{\otimes k}=
\underbrace{V \otimes \ldots \otimes V}_{\text{$k$ times}}.$$
The group $G$ acts in $V^{\otimes k}$ by the $k$-th tensor power of its 
action in $V$: on decomposable tensors we have
$$g(v_1 \otimes \ldots \otimes v_k)=g(v_1) \otimes \ldots \otimes g(v_k).$$
Let us consider the orbit $\bigl\{g v^{\otimes k}: \ g \in G \bigr\}$
for 
$$v^{\otimes k}=
\underbrace{v \otimes \ldots \otimes v}_{\text{$k$ times}}.$$
Then, a linear functional on the orbit of $v^{\otimes k}$ is a polynomial 
of degree $k$ on the orbit of $v$ and hence we can extract some new 
``higher order'' information about the orbit of $v$ by applying already 
developed methods to the orbit of $v^{\otimes k}$. An important observation
is that the orbit $\bigl\{gv^{\otimes k}: \ g \in G \bigr\}$ lies 
in the symmetric part of $V^{\otimes k}$, so the dimension of the affine 
hull of the orbit of $v^{\otimes k}$ does not exceed 
${\dim V + k-1 \choose k}$.

\proclaim{(3.1) Theorem} Let $G$ be a compact group acting in 
a finite-dimensional real vector space $V$, let $v \in V$ be 
a point, and let $\ell: V \longrightarrow {\Bbb R}$ be a linear
functional. Let us define 
$$f: G \longrightarrow {\Bbb R} \qquad \text{by} \qquad 
f(g)=\ell(gv).$$

For an integer $k>0$, let $d_k$ be the dimension of the 
subspace spanned by the orbit $\bigl\{ g v^{\otimes k}: \ g \in G\bigr\}$ 
in $V^{\otimes k}$. In particular, 
$d_k \leq {\dim V + k -1 \choose k}$. Let 
$$\|f\|_{2k}=\left(\int_G f^{2k}(g) \ dg\right)^{1\over 2k}.$$
\roster
\item Suppose that $k$ is odd and that   
$$\int_G f^k(g) \ dg=0.$$
Then 
$$d_k^{-{1 \over 2k}} \|f\|_{2k} \leq \max_{g \in G} f(g) \leq 
d_k^{1 \over 2k} \|f\|_{2k}.$$
\item We have
$$\|f\|_{2k} \leq \max_{g \in G} |f(g)| \leq 
d_k^{1 \over 2k} \|f\|_{2k}.$$  
\endroster
\endproclaim
\demo{Proof} Without loss of generality, we assume that
$f \not\equiv 0$.

Let 
$$B_k(v)=\conv\bigl(gv^{\otimes k}: \ g \in G\bigr)$$
be the convex hull of the orbit of $v^{\otimes k}$. We have 
$\dim B_k(v) \leq d_k$.

Let 
$\ell^{\otimes k} \in \left(V^{\ast}\right)^{\otimes k}$ be the 
$k$-th tensor power of the linear functional $\ell \in V^{\ast}$. 
Thus $f^k(g)=\ell^{\otimes k}\left(g v^{\otimes k}\right)$. 

To prove Part (1), we note that since $k$ is odd,
$$\max_{g \in G} f^k(g)=\left( \max_{g \in G} f(g)\right)^k.$$
Let 
$$u=\int_G g\left(v^{\otimes k}\right) \ dg$$
be the center of $B_k(v)$. Since the average value of $f^k(g)$ is 
equal to 0, we have $\ell^{\otimes k}(u)=0$ and hence
$\ell^{\otimes k}(x)=\ell^{\otimes k}(x-u)$ for all $x \in V^{\otimes k}$.
Let us translate $B_k(v)'=B_k(v)-u$ to the origin and let us consider 
the maximum volume ellipsoid $E$ of the polar of $B_k(v)'$ in its affine 
hull. By Theorem 2.4, we have
$$E=\Bigl\{{\Cal L} \in \left(V^{\otimes k}\right)^{\ast}: \quad 
\int_G {\Cal L}^2\left(gv^{\otimes k}-u\right) \ dg  \leq {1 \over \dim B_k(v)} \Bigr\}.$$

Since the ellipsoid $E$ is contained in the polar 
of $B_k(v)'$, for any 
linear functional ${\Cal L}: V^{\otimes k} \longrightarrow {\Bbb R}$, 
the inequality 
$$\int_G {\Cal L}^2\left(g v^{\otimes k}- u\right) \ dg 
\leq {1 \over d_k} \leq 
{1 \over \dim B_k(v)}$$
implies the inequality
$$ \max_{g \in G} {\Cal L}\left(g v^{\otimes k}-u\right) 
\leq 1.$$
Choosing ${\Cal L}= \lambda \ell^{\otimes k}$ with 
$\lambda=d_k^{-1/2} \|f\|_{2k}^{-k}$, we get the upper bound for
$\max_{g \in G}f(g)$. 

Since the ellipsoid $(\dim E) E$ contains the polar of $B_k(v)'$, 
for any linear functional 
${\Cal L}: V^{\otimes k} \longrightarrow {\Bbb R}$, the inequality
$$ \max_{g \in G} {\Cal L}\left(gv ^{\otimes k}-u\right) \leq 1$$
implies the inequality 
$$\int_G {\Cal L}^2\left(g v^{\otimes k}-u\right) \ dg \leq 
\dim B_k(v) \leq d_k. $$
Choosing ${\Cal L}=\lambda \ell^{\otimes k}$ with any  
$\lambda> \|f\|_{2k}^{-k} d_k^{1/2}$, we obtain the lower bound 
for $\max_{g \in G} f(g)$.

The proof of Part (2) is similar. We modify the definition of 
$B_k(v)$ by letting 
$$B_k(v)=
\conv\bigl(g v^{\otimes k}, -g v^{\otimes k}: \quad g \in G \bigr).$$
The set $B_k(v)$ so defined can be considered as the convex 
hull of an orbit of $G \times {\Bbb Z}_2$ and is centrally symmetric, so  
the ellipsoid $(\sqrt{\dim E})  E$ contains the polar of $B_k(v)$.

Part (2) is also proven by a different method in \cite{Ba02a}.
{\hfill \hfill \hfill} \qed
\enddemo

\remark{Remark} Since $d_k \leq {\dim V + k-1 \choose k}$, the upper 
and lower bounds in Theorem 3.1 are asymptotically equivalent as 
long as $k^{-1} \dim V \longrightarrow 0$. In many interesting cases
we have $d_k \ll {\dim V + k -1 \choose k}$, which results in 
stronger inequalities.
\endremark

\specialhead  Polynomials on the unit sphere \endspecialhead
As is discussed in Examples 1.2 and 2.5, the 
restriction of a homogeneous polynomial 
$f: {\Bbb R}^n \longrightarrow {\Bbb R}$ of degree $m$ onto the 
unit sphere ${\Bbb S}^{n-1} \subset {\Bbb R}^n$ can be viewed as
the restriction of a linear functional 
$\ell: \left({\Bbb R}^n\right)^{\otimes m} \longrightarrow {\Bbb R}$ 
onto the 
orbit of a vector $v=e^{\otimes m}$ for some 
$e \in {\Bbb S}^{n-1}$ in the action of the special 
orthogonal group $SO(n)$. In this case, $v^{\otimes k}=e^{\otimes mk}$ 
spans the symmetric part of $\left({\Bbb R}^n\right)^{mk}$, so 
we have $d_k={n+mk-1 \choose mk}$ in Theorem 3.1. 

Hence Part (1) of Theorem 3.1 implies that if $f$ is an $n$-variate 
homogeneous polynomial of degree $m$ such that 
$$\int_{{\Bbb S}^{n-1}} f^k(x) \ dx=0,$$
where $dx$ is the rotation invariant probability measure on 
${\Bbb S}^{n-1}$, then 
$${n+mk-1 \choose mk}^{-{1 \over 2k}} \|f\|_{2k} 
\leq \max_{x \in {\Bbb S}^{n-1}} f(x) \leq
 {n+mk-1 \choose mk}^{1 \over 2k} \|f\|_{2k},$$
where 
$$\|f\|_{2k}=
\left( \int_{{\Bbb S}^{n-1}} f^{2k}(x) \ dx \right)^{1\over 2k}.$$  
We obtain the following corollary.

\proclaim{(3.2) Corollary} Let us choose $k \geq n \ln(m+1)$. 
Then 
$$\|f\|_{2k} \leq \max_{x \in {\Bbb S}^{n-1}} |f(x)| \leq 
\alpha \|f\|_{2k},$$
for some absolute constant $\alpha>0$ and all homogeneous 
polynomials $f: {\Bbb R}^n \longrightarrow {\Bbb R}$ of degree $m$.
\endproclaim
\demo{Proof}
Applying Part(2) of Theorem 3.1 as above, we 
conclude that for any homogeneous polynomial 
$f: {\Bbb R}^n \longrightarrow {\Bbb R}$ of 
degree $m$,
$$\|f\|_{2k} \leq \max_{x \in {\Bbb S}^{n-1}} |f(x)| \leq 
{n+mk-1 \choose mk}^{1 \over 2k} \|f\|_{2k}$$  
(this inequality is also proven in \cite{Ba02a}).
 
Let
$$H(x)=x \ln {1 \over x} +(1-x) \ln {1 \over 1-x} \quad 
\text{for} \quad 0 \leq x \leq 1$$ 
be the entropy function.
The result now follows from the estimate 
$${a \choose b} \leq \exp\bigl\{a H(b/a) \bigr\},$$
see, for example, Theorem 1.4.5 of \cite{L99}.
{\hfill \hfill \hfill} \qed
\enddemo

Our next application concerns calibrations, see Examples 1.3 and 2.6.
\proclaim{(3.3) Theorem} 
Let $G_m({\Bbb R}^n) \subset \bigwedge^m {\Bbb R}^n$ be the 
Pl\"ucker embedding of the Grassmannian of oriented $m$-subspaces 
of ${\Bbb R}^n$. 
Let $\ell: \bigwedge^m {\Bbb R}^n \longrightarrow {\Bbb  R}$ be 
a linear functional. 
Let 
$$\|\ell\|_{2k}=
\left(\int_{G_m({\Bbb R}^n)} \ell^{2k}(x) \ dx \right)^{1 \over 2k},$$
where $dx$ is the Haar probability measure on $G_m({\Bbb R}^n)$.
Then, for any positive integer $k$,
$$\split &\|\ell\|_{2k} \leq \max_{x \in G_m({\Bbb R}^n)} |\ell(x)| 
\leq (d_k)^{1 \over  2k} \|\ell\|_{2k}, \\
&\text{where} \quad d_k=\prod_{i=1}^m \prod_{j=1}^k 
{n+j-i \over m+k-i-j+1}. \endsplit $$ 
\endproclaim
\demo{Proof} As we discussed in Example 1.3, the Grassmannian 
$G_m({\Bbb R}^n)$ can be viewed as the orbit of 
$v=e_1 \wedge \ldots \wedge e_m$, where $e_1, \ldots, e_n$ is the 
standard basis of ${\Bbb R}^n$, under the action of the 
special orthogonal group $SO(n)$ by the $m$-th exterior power 
of its defining representation in ${\Bbb R}^n$.  
We are going to apply Part (2) of Theorem 3.1 and for that we need to 
estimate the dimension of the subspace spanned by the orbit of
$v^{\otimes k}$. First, we identify $\bigwedge^m {\Bbb R}^n$ with 
the subspace of skew-symmetric tensors in
$\left({\Bbb R}^n\right)^{\otimes m}$ and 
$v$ with the point
$$\sum_{\sigma \in S_m} (\operatorname{sgn\ } \sigma) 
e_{\sigma(1)} \otimes \ldots \otimes e_{\sigma(m)},$$
where $S_m$ is the symmetric group of all permutations of 
$\{1, \ldots, m\}$.

Let us consider $W=({\Bbb R}^n)^{\otimes mk}$. We introduce the 
right action of the symmetric group $S_{mk}$ on $W$ by permutations 
of the factors in the tensor product:
$$\bigl(u_1 \otimes \ldots \otimes u_{mk}\bigr) \sigma=
u_{\sigma(1)} \otimes \ldots \otimes u_{\sigma(mk)}.$$
For $i=1, \ldots, m$, let $R_i \subset S_{mk}$ be the subgroup 
permuting the numbers $1 \leq a \leq mk$ such that $a\equiv i \mod m$ 
and leaving all other numbers intact and for
$j=1, \ldots, k$, let $C_i \subset S_{mk}$ be the subgroup
permuting the numbers $m(i-1)+1 \leq a \leq mi$ and leaving all other 
numbers intact.  

Let $w=e_1 \otimes \ldots \otimes e_m$. 
Then 
$$v^{\otimes k}=(k!)^{-m}w^{\otimes k} 
\left(\sum_{\sigma \in R_1 \times \ldots \times R_m} \sigma \right) 
\left(\sum_{\sigma \in C_1 \times \ldots \times C_k} 
(\operatorname{sgn\ } \sigma) \sigma \right).$$
It follows then that $v^{\otimes k}$ generates the $GL_n$-module 
indexed by the rectangular $m \times k$ Young diagram, so its dimension
$d_k$ is given by the formula of the Theorem, see Chapter 6 of
\cite{FH91}.    
{\hfill \hfill \hfill} \qed 
\enddemo

\proclaim{(3.4) Corollary} Under the conditions of Theorem 3.3,
let $k \geq mn \ln (m+1)$. Then
$$\|\ell\|_{2k} \leq \text{comass of\ } \ell \leq \alpha \|\ell\|_{2k}$$
for some absolute constant $\alpha>0$.
\endproclaim
\demo{Proof} 
We have 
$$\split d_k &\leq  \prod_{i=1}^m \prod_{j=1}^k {n+j-i \over k-j+1} 
\leq \left(\prod_{j=1}^k {n+j-1 \over k-j+1}\right)^m=
{n+k-1 \choose n-1}^m \\
&\leq \exp\left\{m(n+k-1) H\left( {n-1 \over n+k-1} \right) \right\},
\endsplit$$
cf. Corollary 3.2. The proof now follows.
{\hfill \hfill \hfill} \qed
\enddemo

To understand the convex geometry of an orbit, we would like to compute 
the maximum value of a ``typical'' linear functional on the 
orbit. Theorem 3.1 allows us to replace the maximum value by an 
$L^p$ norm. To estimate the average value of 
an $L^p$ norm, we use the following simple computation.
\proclaim{(3.5) Lemma} Let $G$ be a compact group acting in a 
$d$-dimensional real vector space $V$ endowed with a $G$-invariant 
scalar product $\langle \cdot, \cdot \rangle$ and let 
$v \in V$ be a point. Let 
${\Bbb S}^{d-1} \subset V$ be the unit sphere endowed with the 
Haar probability measure $dc$. Then, for every positive integer $k$,
we have
$$\int_{{\Bbb S}^{d-1}} \left(\int_G \langle c, gv \rangle^{2k}\ dg
\right)^{1 \over 2k} \ dc \leq  
\sqrt{{2k \langle v, v \rangle \over d}}
.$$
\endproclaim
\demo{Proof} Applying H\"older's inequality, we get 
$$\int_{{\Bbb S}^{d-1}} \left(\int_G \langle c, gv \rangle^{2k}\ dg
\right)^{1 \over 2k} \ dc \leq  
\left(\int_{{\Bbb S}^{d-1}} \int_G \langle c, gv \rangle^{2k} \ dg \ dc
\right)^{1 \over 2k}.$$
Interchanging the integrals, we get  
$$\int_{{\Bbb S}^{d-1}} \int_G \langle c, gv \rangle^{2k} \ dg \ dc=
\int_G \left(\int_{{\Bbb S}^{d-1}} \langle c, gv \rangle^{2k} \ dc
\right) \ dg. \tag3.5.1$$
Now we observe that the integral inside has the same value for 
all $g \in G$. Therefore, (3.5.1) is equal to 
$$\int_{{\Bbb S}^{d-1}} \langle c, v \rangle^{2k} \ dc
 =\langle v, v \rangle^k 
{\Gamma(d/2) \Gamma(k+1/2) \over \sqrt{\pi} \Gamma(k+ d/2)},$$
see, for example, \cite{Ba02a}.  

Now we use that $\Gamma(k+1/2) \leq \Gamma(k+1) \leq k^k$ 
and 
$${\Gamma(d/2) \over \Gamma(k+d/2)}=
{1 \over (d/2) (d/2+1) \cdots (d/2+k-1)} \leq (d/2)^{-k}.$$ 
{\hfill \hfill \hfill } \qed 
\enddemo

\head 4. Some Geometric Corollaries \endhead

\specialhead The metric structure of the unit comass ball 
\endspecialhead
Let $V_{m,n}=\bigwedge^m {\Bbb R}^n$ with the orthonormal basis 
$e_I=e_{i_1} \wedge \ldots \wedge e_{i_m}$, where 
$I$ is an $m$-subset $1 \leq i_1 < i_2 < \ldots < i_m \leq n$ of 
the set $\{1, \ldots, n\}$, and the corresponding scalar 
product $\langle \cdot, \cdot \rangle$.
Let $G_m({\Bbb R}^n) \subset V_{m,n}$ be the Pl\"ucker embedding of the 
Grassmannian of oriented $m$-subspaces of ${\Bbb R}^n$,
let $B=\conv\left(G_m({\Bbb R}^n)\right)$ be the unit mass ball, 
and let $B^{\circ}  \subset V_{m,n}^{\ast}=V_{m,n}$ be the 
unit comass ball, consisting of the linear functionals 
with the maximum value on $G_m({\Bbb R}^n)$ not exceeding 1, 
see Examples 1.3 and 2.6.

The most well-known example of a linear functional 
$\ell: V_{m,n} \longrightarrow {\Bbb R}$ of comass 1 is given 
by an exterior power of the K\"ahler form. Let us suppose that
$m$ and $n$ are even, so $m=2p$ and $n=2q$. Let 
$$\split
 &\omega=e_1 \wedge e_2 + e_3 \wedge e_4 + \ldots + e_{q-1} \wedge e_q
\quad \text{and} \\ 
 &f={1 \over p!}
\underbrace{\omega \wedge \ldots \wedge \omega}_{\text{$p$ times}} 
\in V_{m,n}. \endsplit$$
Then 
$$\max_{x \in G_m({\Bbb R}^n)} \langle f, x \rangle=1,$$
and, moreover, the subspaces $x \in G_m({\Bbb R}^n)$ where the 
maximum value 1 is attained look as follows. We identify 
${\Bbb R}^n$ with ${\Bbb C}^q$ by identifying 
$${\Bbb R}e_1 \oplus {\Bbb R}e_2={\Bbb R}e_3 \oplus {\Bbb R}e_4 
=\ldots ={\Bbb R} e_{q-1} \oplus {\Bbb R} e_q ={\Bbb C}.$$
Then the subspaces $x \in G_m({\Bbb R}^n)$ with 
$\langle f,x \rangle=1$ are exactly those identified with the complex 
$p$-dimensional subspaces of ${\Bbb C}^q$, see \cite{HL82}. 

We note that the Euclidean length 
$\langle f, f \rangle^{1/2}$ of $f$ is equal to 
$\displaystyle {q \choose p}^{1/2}$. 
In particular, if $m=2p$ is fixed and $n=2q$
grows, the length of $f$ grows as $n^{p/2}=n^{m/4}$. 

Another example is provided by the special Lagrangian calibration $a$.
In this case, $n=2m$ and 
$$a=\operatorname{Re\ } (e_1+ie_2) \wedge \ldots 
\wedge (e_{2m-1} +i e_{2m}).$$
The length $\langle a, a \rangle^{1/2}$ of $a$ is 
$\displaystyle \left(\sum_{j \leq m/2} {m \choose 2j}\right)^{1/2}$. 
The maximum value of $\langle a, x \rangle$ for $x \in G_m({\Bbb R}^n)$ is 
1 and it is attained on the ``special Lagrangian subspaces'', see 
\cite{HL82}.       
 
The following result shows that there exist calibrations with a much 
larger Euclidean length than that of the power $f$ of the K\"ahler form 
or the special Lagrangian calibration $a$.
\proclaim{(4.1) Theorem} 
\roster
\item Let $c \in V_{m,n}$ be a vector such that
$$\max_{x \in G_m({\Bbb R}^n)} \langle c, x \rangle =1.$$
Then 
$$\langle c, c \rangle \leq {n \choose m}.$$
\item There exists $c \in V_{m,n}$ such that 
$$\max_{x \in G_m({\Bbb R}^n)} \langle c, x \rangle=1$$
and
$$\langle c, c \rangle \geq \alpha 
\bigl(nm \ln(m+1) \bigr)^{-1} {n \choose m},$$
where $\alpha>0$ is an absolute constant.
\endroster
\endproclaim
\demo{Proof}
Part (1) follows since the convex hull of the Grassmannian contains 
a ball of radius $\displaystyle {n \choose m}^{-1/2}$, see Example 2.6.

To prove Part (2), let us choose $k=\lceil mn\ln(m+1) \rceil$ in 
Lemma 3.5. Then, by Corollary 3.4, the maximum value of 
$\langle c, x \rangle$ for $x \in G_m({\Bbb R}^n)$ is approximated 
by 
$$\left(\int_{G_m({\Bbb R}^n)} 
\langle c, x \rangle^{2k} \ dx \right)^{1\over 2k}$$ within 
a constant factor. We apply Lemma 3.5 with $V=V_{m,n}$,
$d={n \choose m}$, $G=SO(n)$, and $v=e_1 \wedge \ldots \wedge e_m$.
Hence $\langle v, v \rangle=1$ and there exists $c \in V_{m,n}$ 
with $\langle c, c \rangle=1$ and such that 
$$\left(\int_{G_m({\Bbb R}^n)} 
\langle c, x \rangle^{2k} \ dx \right)^{1\over 2k} \leq 
\sqrt{2k} {n \choose m}^{-1/2}.$$
Rescaling $c$ to a comass 1 functional, we 
complete the proof of Part (2). 
{\hfill \hfill \hfill} \qed
\enddemo

For $m=2$ the estimate of Part (2) is exact up to an absolute 
constant, as witnessed by the K\"ahler calibration. However, for 
$m \geq 3$, the calibration $c$ of Part (2) has a larger length 
than the K\"ahler or special Lagrangian calibrations. The gap 
only increases when $m$ and $n$ grow. The distance to the origin of the 
supporting hyperplane $\langle c, x \rangle=1$ of the face 
of the convex hull of the Grassmannian is 
equal to $\langle c, c \rangle^{-1/2}$ so the faces defined by 
longer calibrations $c$ are closer to the origin. Thus, the 
faces spanned by complex subspaces or the faces spanned by 
special Lagrangian subspaces are much more ``shallow'' than the 
faces defined by calibrations $c$ in Part (2) of the Theorem.
We do not know if those ``deep'' faces are related to any interesting 
geometry. Intuitively, the closer the face to the origin, the 
larger piece of the Grassmannian it contains, so it is 
quite possible that some interesting classes of manifolds are 
associated with the ``long'' calibrations $c$. 
  
\specialhead The volume of the set of non-negative polynomials 
\endspecialhead

Let $U_{m,n}$ be the space of real homogeneous polynomials $p$ of 
degree $m$ in $n$ variables such that the average value of 
$p$ on the unit sphere ${\Bbb S}^{n-1} \subset {\Bbb R}^n$ is 0, 
so $\dim U_{m,n} ={n+m-1 \choose m}-1$ for $m$ even and 
$\dim U_{m,n}={n+m-1 \choose m}$ for $m$ odd. As before, we make 
$U_{m,n}$ a Euclidean space with the $L^2$ inner product 
$$\langle f, g \rangle=\int_{{\Bbb S}^{n-1}} f(x) g(x) \ dx.$$
We obtain the following corollary.
\proclaim{(4.2) Corollary} Let $\Sigma_{m, n} \subset U_{m,n}$ 
be the unit sphere, consisting of the polynomials 
with $L^2$-norm equal to 1. For a polynomial $p \in U_{m,n}$, let
$$\|p\|_{\infty}=\max_{x \in {\Bbb S}^{n-1}} |p(x)|.$$
Then 
$$\int_{\Sigma_{m,n}} \|p\|_{\infty} \ dp \leq \alpha \sqrt{n \ln(m+1)}$$
for some absolute constant $\alpha>0$.
\endproclaim 
\demo{Proof} Let us choose $k=\lceil n\ln(m+1) \rceil$. Then, by 
Corollary 3.2, $\|p\|_{2k}$ approximates 
$\|p\|_{\infty}$ within an absolute constant.  

Now we use Lemma 3.5. As in Examples 1.2 and 2.5, we identify space
$U_{m,n}$ with the space of linear functionals $\langle c, gv \rangle$ 
on the orbit $\bigl\{gv: \ g \in SO(n) \bigr\}$ of 
$v$. By the remark after the proof of Theorem 2.4, we have 
$\langle v, v \rangle =\dim U_{m,n}$. The proof now follows. 
{\hfill \hfill \hfill} \qed
\enddemo

Thus the $L^{\infty}$-norm of a typical $n$-variate polynomial of 
degree $m$ of the unit $L^2$-norm in $U_{m,n}$ is 
$O\bigl(\sqrt{n \ln(m+1)}\bigr)$.
In contrast, the $L^{\infty}$ norm of a {\it particular} polynomial 
can be of the order of $n^{m/2}$, that is, substantially bigger.

Corollary 4.2 was used by the second author to obtain a bound on the 
volume of the set of non-negative polynomials. 

Let us consider the shifted set $Pos_{2k,n}' \subset U_{2k,n}$
of non-negative polynomials defined by (1.2.2). 
We measure the size of a set $X \subset U_{2k,n}$ by the quantity
$\displaystyle \left({\vl X \over \vl K } \right)^{1/d}$, 
where $d=\dim U_{2k,n}$ and 
$K$ is the unit ball in $U_{2k,n}$, which is more ``robust'' than just
the volume $\vl X$, as it takes into account the effect of a high
dimension, cf. Chapter 6 of \cite{P89}. 

The following result is from \cite{Bl03}, we made some trivial 
improvement in the dependence on the degree $2k$.
\proclaim{(4.3) Theorem} Let $Pos_{2k,n}' \subset U_{2k,n}$ be the
shifted set of non-negative polynomials, let $K \subset U_{2k,n}$ be the 
unit ball and let $d=\dim U_{2k,n}={n + 2k -1 \choose 2k}-1$. Then 
$$\left( \vl Pos_{2k,n} \over \vl K \right)^{1/d} \geq 
\alpha \bigl(n \ln(2k+1) \bigr)^{-1/2}$$
for some absolute constant $\alpha>0$. 
\endproclaim 
\demo{Proof} Let $\Sigma_{2k,n} \subset U_{2k,n}$ be the unit sphere.
Let $p \in \Sigma_{2k,n}$ be a point. The ray $\lambda p: \lambda \geq 0$ 
intersects the boundary of $Pos_{2k,n}'$ at a point $p_1$ such that 
$\min_{x \in {\Bbb S}^{n-1}} p_1(x)=-1$, so the length of the 
interval $[0, p_1]$ is 
$|\min_{x \in {\Bbb S}^{n-1}} p(x)| \leq \|p\|_{\infty}$.

Hence 
$$\split &\left( \vl Pos_{2k,n}' \over \vl K \right)^{1/d} =
\left(\int_{\Sigma_{2k,n}} |\min_{x \in {\Bbb S}^{n-1}} p(x)|^{-d} \ dp
\right)^{1/d} \geq  
\left(\int_{\Sigma_{2k,n}} \|p\|^{-d}_{\infty} \ dp
\right)^{1/d} \\ \geq 
&\int_{\Sigma_{2k,n}} \|p\|_{\infty}^{-1} \ dp 
\geq \left(\int_{\Sigma_{2k,n}} \|p\|_{\infty} \ dp
\right)^{-1}, \endsplit $$
by the consecutive application of H\"older's and Jensen's inequalities,
so the proof follows by Corollary 4.2.
{\hfill \hfill \hfill} \qed
\enddemo   

We defined $Pos_{2k,n}$ as the 
set of non-negative polynomials with the average value 1 on the 
unit sphere, see (1.2.1). 
There is an important subset $Sq_{2k, n} \subset Pos_{2k,n}$, 
consisting of the polynomials that are sums of squares of homogeneous 
polynomials of degree $k$. It is known that $Pos_{2k,n}=Sq_{2k,n}$ 
if $k=1$, $n=2$, or $k=2$ and $n=3$, see Chapter 6 of \cite{BCR98}.
The following result from \cite{Bl03}
shows that, in general, $Sq_{2k, n}$ is a rather small subset 
of $Pos_{2k,n}$.
 
Translating $p \longmapsto p- (x_1^2 + \ldots +x_n^2)^k$, 
we identify $Sq_{2k,n}$ with a subset $Sq_{2k,n}'$ of $U_{2k, n}$.
\proclaim{(4.4) Theorem} Let $Sq_{2k,n}' \subset U_{2k,n}$ be the 
shifted set of sums of squares, let $K \subset U_{2k,n}$ be the 
unit ball and let $d=\dim U_{2k,n}={n + 2k -1 \choose 2k}-1$. Then 
$$\left( \vl Sq_{2k,n}  \over \vl K  \right)^{1/d} \leq 
\alpha 2^{4k} {n+k-1 \choose k}^{1/2} 
{n+2k-1 \choose 2k}^{-1/2}$$
for some absolute constant $\alpha>0$.
\endproclaim 

In particular, if $k$ is fixed and $n$ grows, the upper bound has the 
form $c(k) n^{-k/2}$ for some $c(k)>0$.

The proof is based on bounding the right hand side 
of the inequality of Theorem 4.4 by the average 
width of $Sq_{2k,n}$, cf. Section 6.2 of \cite{S93}. The average 
width is represented by the integral
$$\int_{\Sigma_{2k, n}}  
\left(\max_{f \in \Sigma_{k,n}} \langle g, f^2 \rangle - 
\min_{f \in \Sigma_{k,n}} \langle g, f^2 \rangle \right) \ dg.$$
By Corollary 3.2, we can bound the integrand by 
$$\alpha \left(\int_{\Sigma_{k,n}} \langle g, f^2 \rangle^{2q} 
\ df\right)^{1 \over 2q}$$
for some absolute constant $\alpha$ and $q={n+k-1 \choose k}$ 
and proceed as in the proof of Lemma 3.5. 

\head Acknowledgment \endhead

We thank B. Sturmfels for suggesting to us to consider the convex hull 
of the Grassmannian and J. Sullivan for pointing to connections 
with calibrated geometries. 

\refstyle{A}

\enddocument

\end